\documentclass[11pt%
							, draft%
								]{amsart}
\usepackage{research}

\usepackage{fullpage,color,url}
\usepackage{pdfsync}

\usepackage[square,comma,sort&compress]{natbib}
\setcitestyle{numbers}
\usepackage{ifpdf}
\ifpdf
\usepackage{hyperref}
\fi

\numberwithin{equation}{section}
\numberwithin{figure}{section}

\nc{\Gz}{G(\{0\})}
\nc{\eb}{y}
\nc{\ebt}{\eb(t)}

\newcommand{\ecn}[1]{#1^{(b_n)}}
\nc{\tat}{\tau(t)}
\nc{\tbt}{\tau(b(t))}
\nc{\tbn}{\tau(b_n)}

\nc{\Gg}{G_\gamma}
\nc{\Eg}{\E_{\gamma}}
\nc{\Egc}{\cl{\E}_{\gamma}}

\newcommand{\rinv}[1]{#1^{-\rho}}

\nc{\limpp}{\zeta}
\nc{\nua}{\nu_\alpha}
\nc{\tk}{\tilde{\tau}_{k+1}(b_n)}
\nc{\tkit}{\tilde{\tau}_k(t)}
\nc{\tkoit}{\tilde{\tau}_{k+1}(t)}
\nc{\tsk}{\tau^*_{k+1}}
\nc{\sk}{\sigma_{k+1}(b_n)}
\nc{\skt}{\sigma_{k+1}(t)}
\nc{\ssk}{\sigma^*_{k+1}}
\nc{\Xkj}{X_{kj}}
\nc{\eXkj}{\ecn{\Xkj}}
\nc{\Xko}{X_{k0}}
\nc{\Xkm}{X_{km}}
\nc{\eXkm}{\ecn{X_{km}}}
\nc{\Xvkm}{\ecn{\bX_{k,m}}}
\nc{\xikj}{\xi_{k}(j)}
\nc{\xiko}{\xi_{k}(0)}
\nc{\xikm}{\xi_{k}(m)}
\nc{\xivkm}{\bxi_{k,m}}
\nc{\ppnveco}{\vartheta_n}
\nc{\ppveco}{\vartheta}
\nc{\ppnvec}{\eta'_n}
\nc{\ppvec}{\eta'}
\nc{\Ld}{\Lambda_\delta}
\nc{\ppnvst}{\hat{\eta}_n}
\nc{\ppvst}{\hat{\eta}}
\nc{\ppnvrest}{\eta''_n}
\nc{\ppvrest}{\eta''}
\nc{\ppnst}{\tilde{\eta}_n}
\nc{\ppst}{\tilde{\eta}}
\nc{\ppnstj}{\tilde{\eta}^*_n}
\nc{\ppnmain}{\eta_n}
\nc{\ppmain}{\eta}
\nc{\ppnfull}{\eta^*_n}
\nc{\ppfull}{\eta^*}

\nc{\ta}{\tau_A}
\nc{\taa}{\tau^A}
\nc{\taz}{\taa_0}
\nc{\tak}{\taa_k}
\nc{\tako}{\taa_{k+1}}
\nc{\tkt}{\tau_k(t)} 
\nc{\tkot}{\tau_{k+1}(t)} 
\nc{\tkob}{\tau_{k+1}(b_n)}
\nc{\tkb}{\tau_k(b_n)}
\nc{\Kc}{K'}
\nc{\Xc}{X'}
\nc{\bXc}{\bX'}
\nc{\tac}{\tau'_A}
\nc{\ttc}{\tau'(t)}
\nc{\tbtc}{\tau'(b(t))}

\nc{\ppmcfull}{\chi^{**}_n}
\nc{\ppmco}{\chi_n^0}
\nc{\ppmc}{\chi^*_n}
\nc{\ppmcd}{\chi_n}
\nc{\Nn}{N_n}
\nc{\Nnd}{\tilde{N}_n}
\nc{\ppmcde}{\chi'_n}
\nc{\ppmcdc}{\chi''_n}
\nc{\ppmce}{{\chi^*_n}'}
\nc{\ppmcc}{{\chi^*_n}''}

\def\binfty{\boldsymbol \infty}

\title[Markov Kernels and CEV]
{Markov Kernels and the Conditional Extreme Value Model}
\author[S.~I.~Resnick]{Sidney I.~Resnick}
\author[D.~Zeber]{David Zeber}

\address{Sidney I.~Resnick, School of ORIE, Rhodes Hall 284, Cornell University,  Ithaca, NY 14853}
\email{sir1@cornell.edu}%
\address{David Zeber, Department of Statistical Science, 301 Malott Hall, Cornell University, Ithaca NY 14853}
\email{dsz5@cornell.edu}%
 
\thanks{S.~I.~Resnick and D.~Zeber were partially supported by ARO   Contract W911NF-10-1-0289 and NSA Grant H98230-11-1-0193 at Cornell University.}

\date{\today}

\begin{document}

\begin{abstract}
  The classical approach to multivariate extreme value modelling
  assumes
 that the joint distribution belongs to a
  multivariate domain of attraction. This requires 
  each marginal distribution be individually attracted to a univariate
  extreme value distribution.  
An apparently more flexible extremal model for multivariate data
   was proposed by Heffernan and Tawn
   under which not all the components are required to belong to
  an extremal domain of attraction but assumes instead the
  existence of an asymptotic approximation to the conditional
  distribution of the random vector given one of the components
is extreme. Combined with the knowledge that the conditioning
  component belongs to a univariate domain of attraction, this leads
  to an approximation of the probability of certain risk regions.
The original focus on conditional distributions had
technical
drawbacks
 but is natural in several contexts.  We place this approach in the
 context of the more
general approach using convergence of measures and multivariate
regular variation on cones.\end{abstract}
%


\maketitle


\section{Overview}

The classical approach to extreme value modelling for multivariate data is to assume that the joint distribution belongs to a multivariate domain of attraction. 
In particular, this requires that each marginal distribution be individually attracted to a univariate extreme value distribution. 
The domain of attraction condition may be phrased conveniently in terms of regular variation of the joint distribution on an appropriate cone; see Das and Resnick \cite[Proposition 4.1]{das2011conditioning}.

A more flexible model for data realizations of a random vector was proposed by Heffernan and Tawn \cite{heffernan2004conditional}, under which not all the components are required to belong to an extremal domain of attraction. 
Such a model accomodates varying degrees of asymptotic dependence between pairs of components. 
Instead of starting from the joint distribution, Heffernan and Tawn
assumed the existence of an asymptotic approximation to the
conditional distribution of the random vector given one of the
components was extreme.  
Combined with the knowledge that the conditioning component belongs to
a univariate domain of attraction, this leads to an approximation for
the probabilities of certain multivariate risk sets.
However,  focusing  on conditional distributions creates problems when
taking limits owing to ambiguity regarding the choice of version.  So 
the Heffernan/Tawn approach was reformulated as the
{\it conditional extreme value model\/} (CEVM) in 
\cite{heffernan2007limit,
das2011conditioning, das2011detecting} using regular
variation of the joint distributions on a smaller
cone than the one employed in multivariate extreme value theory, an
approach related to hidden regular variation
\cite{resnick2008multivariate,maulik2004characterizations, 
resnick2002hidden,resnick2007heavy}.  

Conditional distributions are natural objects in many
circumstances,
for example if densities exist or if one variable is an explicit
function of others. So
we return to the Heffernan and Tawn
\cite{heffernan2004conditional} formulation,
placing it in a formal context that uses the idea of 
transition kernels in a domain of attraction developed in
\cite{resnick2011asymptotics,resnick:zeber:2013b}. We see how reliance
on transition kernels fits with general theory expressed in terms of
vague convergence of measures and to what extent the reliance on
kernels restricts the class of limit measures. 

In order to better fit in with the study of extremes of a
random vector, we extend the kernel domain of attraction condition
used in \cite{resnick2011asymptotics} beyond standardized regular
variation
 to accomodate general linear normalization in both
the initial state and the distribution of the next state.  We examine
conditions under which this extends to a CEVM, when combined with a
marginal domain of attraction assumption, and we derive explicit
formulas for the CEV limit measure in different cases.  Also, through
a number of revealing examples, we explore the properties of the
normalization functions, and technicalities surrounding the choice of
version of the conditional distribution and the limit distribution $G$.

Section \ref{sec:back} summarizes necessary background, definitions
and basic results including when a random vector $(X,Y)$ satisfies the
{\it conditional extreme value model\/} (CEVM). Section
\ref{stdcontbasiccase} treats the {\it standard\/} case where both $X$
and $Y$ can be scaled by the same function and this restriction is
weakened in Sections \ref{sectiongensc} and \ref{seccevmgeny}.

We denote by $ \mplus(\E)$ the space of Radon measures on a nice space
$\E$ topologized by vague convergence which is written as
$\stackrel{v}{\to}$ \cite{resnick2007extreme}. 
Weak convergence  \cite{billingsley1999convergence} 
of probability measures is
denoted $\Rightarrow$.
We write $\xi \sim G$ to mean a random variable
$\xi$ has distribution $G(\cdot) $ and, if no confusion can arise, we
often use $G$ to also mean 
the distribution function $\P[\xi \leq x]=G(x)$. Alternatively for a
random variable $Y$, we write $F_Y$ for the distribution of $Y$.
The class of regularly varying functions with index $\rho$ on
$(0,\infty)$ is $RV_\rho$
\cite{bingham1989regular,feller:1971,de2006extreme,resnick2007extreme}. The
probability measure degenerate at $c\in \mathbb{R}$ is $\epsilon_c(\cdot).$


\section{Background}\label{sec:back}

First, we review the basics of 
extended regular variation, which features prominently in the
formulation of the CEVM, as well as
concepts of univariate
extreme value theory.  
We then define the conditional extreme value model  and discuss its basic properties.

\subsection{Extended Regular Variation} 
\label{secerv}

Regular variation and extended regular variation is important
 in the mathematical description of extreme  and conditional
 extreme value theory
\cite{bingham1989regular, resnick2007extreme, resnick2007heavy, seneta1976regularly,
de2006extreme}. 
The pair of functions $a : (0,\infty) \mapsto (0,\infty)$ and $f :
(0,\infty) \mapsto \R$ are \emph{extended regularly varying} (ERV)
with parameters $\rho, k \in \R$ if as $t\to\infty$,
\begin{equation} \label{eqdeferv}
\frac{a(tx)}{a(t)}\conv x^{\rho} \qquad\text{and}\qquad \frac{f(tx)-f(t)}{a(t)}\conv\psi(x), \qquad x>0, 
\end{equation} \cite[Appendix B.2]{de2006extreme},
where 
\begin{equation} \label{eqervpsi}
\psi(x) = \begin{cases} k\inv{\rho}(x^{\rho}-1) & \rho \neq 0 \\ k \log x & \rho=0 \end{cases}. 
\end{equation} 
We will write this as $a,f \in \ERV_{\rho,k}$ with  $a \in \RV_\rho$.
A useful identity is 
\begin{equation} \label{eqervpsiinv}
\psi(\inv{x}) = -\rinv{x} \psi(x). 
\end{equation}
Note that this differs slightly from the usual definition of extended regular variation, which assumes $k=1$. 
If $\phi(x) := \limt (f(tx)-f(t))/a(t)$ exists for $x>0$, then $a$ is necessarily regularly varying, and $\phi \equiv \psi$, the function given in \eqref{eqervpsi}. Also, the convergences in \eqref{eqdeferv} are locally uniform, implying that 
\begin{equation*} 
\frac{a(tx_t)}{a(t)}\conv x^{\rho} \qquad\text{and}\qquad \frac{f(tx_t)-f(t)}{a(t)}\conv\psi(x) \mt{whenever} x_t \goesto x>0.  
\end{equation*}

Furthermore, if $k\neq 0$ we obtain the following properties depending on the value of $\rho$. Recall the sign function $\sgn(u) = u/\abs{u} \ind{u \neq 0}$. 
\begin{itemize}
	\item  If $\rho > 0$, then $f \cdot \sgn(k) \in \RV_\rho$, and $f(t)/a(t) \goesto k/\rho$. 
	\item If $\rho < 0$, then $f(\infty)= \lim_{t\goesto\infty}
          f(t)$ exists finite, $(f(\infty)-f(t))/a(t) \goesto
          k/\abs{\rho}$ and   $(f(\infty)-f) \cdot \sgn(k) \in
          \RV_{-\abs{\rho}}$.
	\item If $\rho = 0$, \ie $a$ is slowly varying, then $f \in
          \Pi(a)$,
the $\Pi$-varying functions with auxiliary function $a(\cdot)$
\cite[Appendix B.2]{de2006extreme}). Suppose $k>0$. Then $f(\infty)
\leq \infty$ exists. If $f(\infty) = \infty$, then $f\in \RV_0$ and
$f(t)/a(t) \goesto\infty$. If $f(\infty) < \infty$, then $f(\infty)-f
\in \RV_0$, and ${(f(\infty)-f(t))/a(t)} \goesto \infty$. If $k<0$,
then $-f$ has these properties.  
\end{itemize}


\subsection{Domains of Attraction}

\label{secdoa}

 For $\gamma \in \R$, define $\Eg = \{x \in \R : 1 + \gamma x > 0\}$
 so that
\begin{equation} \label{eqegamma}
 \Eg = \begin{cases} (-\inv{\gamma},\infty) & \gamma > 0 \\ (-\infty,\infty) & \gamma = 0 \\ (-\infty, \inv{\abs{\gamma}}) & \gamma < 0 \end{cases}. 
\end{equation}
The distribution $F_Y$ of a random variable $Y$ is in the \emph{domain of attraction} of an extreme value distribution $\Gg$ for some $\gamma\in\R$, written $F_Y\in D(\Gg)$, if there exist functions $a(t)>0$ and $b(t) \in \R$ such that 
\[ F_Y^t\big(a(t)y + b(t)\big) \conv \Gg(y) \] 
weakly as $t\goesto\infty$, where $\Gg(y) = \exp\{-(1+\gamma
y)^{-1/\gamma}\}$ for $y\in \Eg$ \cite{de2006extreme,resnick2007extreme}.
This can be reformulated in terms of the tail of the distribution $F_Y$ as 
\begin{equation}\label{eqdoagen}
t\P\bigl[{Y-b(t)}/{a(t)}>y \bigr] \conv (1+\gamma y)^{-1/\gamma}, \qquad y\in\Eg.
\end{equation}
If $\gamma=0$, we interpret the limit as $e^{-y}$. 

If \eqref{eqdoagen} holds for some functions $a$ and $b$, then it
holds for 
\cite[Theorem 1.1.6,
p.~10]{de2006extreme}) 
\begin{equation} \label{eqbtdf}
 b(t) = \ginv{\bigg(\frac{1}{1-F_Y} \bigg)}(t) = \ginv{F}_Y\big(1-\inv{t}\big), 
\end{equation}
where $\ginv{g}$ is the left-continuous inverse of the nondecreasing function $g$. 
By inversion, \eqref{eqdoagen} yields
\begin{equation} \label{eqbterv}
 \frac{b(tx)-b(t)}{a(t)} \goesto \frac{x^\gamma-1}{\gamma} \ind{\gamma \neq 0} + \log x \ind{\gamma = 0}, 
\end{equation}
\ie $a,b \in \ERV_{\gamma,1}$. 
Furthermore, if functions $\tilde{a}>0$ and $\tilde{b}\in \R$ on $(0,\infty)$ are \emph{asymptotically equivalent} to $a,b$, \ie they satisfy 
\[ \frac{\tilde{a}(t)}{a(t)} \conv 1 \mand \frac{{\tilde{b}}(t)-b(t)}{a(t)} \conv 0 \mt{as} t\goesto\infty, 
\]
then \eqref{eqdoagen} and \eqref{eqbterv} hold with $a,b$ replaced by $\tilde{a},\tilde{b}$. It follows that  \eqref{eqdoagen} is equivalent to 
$t \P[ \ginv{b}(Y) > ty ] \goesto \inv{y}$ for $y>0$, 
\ie $1-F_{\ginv{b}(Y)} \in \RV_{-1}$. This is known as standardization (see \cite[Chapter 5]{resnick2007extreme}). 
We say that $Y^*$ is in the \emph{standardized domain of attraction}, and write $F_{Y^*} \in D(G^*_1)$, if 
\begin{equation}\label{eq:tailconv}
 t \P[ Y> ty] \conv \inv{y}, \qquad y>0, \end{equation}
a variant of \eqref{eqdoagen} for $\gamma=1$.

\subsection{The Conditional Extreme Value (CEV) Model}

Denote by $\Egc$ the closure on the right of the interval $\Eg$. 
A bivariate random vector $(X,Y)$ on $\R^2$ follows a
\emph{conditional extreme value model} (CEVM) if there exists a
 measure $\mu \in \mplus ([-\infty,\infty]\times \Egc)$, and
functions $a(t),\alpha(t)>0$, $b(t),\beta(t) \in \R$, such that, as
$t\goesto\infty$, 
\begin{equation} \label{eqcevmdef}
 t\P\bigg[ \bigg(\frac{X-\beta(t)}{\alpha(t)}, \frac{Y-b(t)}{a(t)} \bigg) \in \cdot\, \bigg] \vconv \mu(\cdot) 
\mt{in} \mplus([-\infty,\infty]\times \Egc),
\end{equation}
and where $\mu$ satisfies the \emph{conditional non-degeneracy} conditions: for each $y \in \Eg$, 
\begin{equation} \label{eqcondnondegennew}
\begin{split}
&\mu([-\infty,x] \times (y,\infty] ) \text{ is not a degenerate distribution in } x ; \\
& \mu(\{\infty\} \times (y,\infty]) = 0.  
\end{split}
\end{equation}
It is convenient to choose the normalization such that 
\begin{equation} \label{eqcevdisth}
H(x) := \mu([-\infty,x] \times (0,\infty] ) \text{ is a probability distribution on } [-\infty,\infty]. 
\end{equation}
See \cite{heffernan2007limit,das2011conditioning} for details and
\cite{heffernan2004conditional} for background.

Some remarks: By applying the joint convergence \eqref{eqcevmdef} to
rectangles $[-\infty,\infty]\times (y,\infty]$, we see that the
distribution of $Y$ is necessarily attracted to $\Gg$ for some
$\gamma$.  Also, an important property is that the functions
$\alpha,\beta$ are ERV for some $\rho,k \in \R$
\mbox{\cite[Proposition 1]{heffernan2007limit}}.  The limit measure
$\mu$ in \eqref{eqcevmdef} is a product measure if and only if
$(\rho,k) = (0,0)$ \cite[Proposition 2]{heffernan2007limit}.

Condition \eqref{eqcondnondegennew}
is somewhat different from  what is given in 
\cite{das2011conditioning,heffernan2007limit} which 
failed to preclude mass on the line $\{\infty\}\times (-\infty,
\infty]$
through infinity. Mass on this line invalidates the convergence to
types theorem \cite{resnick:1999book, feller:1971} and since the
theory in \cite{das2011conditioning,heffernan2007limit} employs 
convergence of types arguments, we require the second condition in
\eqref{eqcondnondegennew}.
Condition \eqref{eqcevmdef} entails 
$Y\in D(G_\gamma)$ and  $\mu([-\infty,x] \times \{\infty\}) =
0$.
 Example \ref{exdefecg}
is a case where \eqref{eqcevmdef} holds for two distinct
 normalizations, which are not asymptotically equivalent,
 yielding two distince limit measures.
 One limit measures has $\mu(\{\infty\} \times (y,\infty]) > 0$
 and the other has $\mu(\{\infty\} \times (y,\infty]) =0.$


\section{Standard Case}
\label{stdcontbasiccase}

Let $(X,Y)$ be a random vector on $\R^2$, with dependence specified by a transition kernel $K$: 
\[ \P[X \in A \bgiven Y = y ] = \mckern{y}{\cdot} \qquad y \in \R. \]
$K(y,A)$ is a measure in the second variable $A$ and measurable in $y$
for each fixed $A$.
We show if the distribution of $Y$ is in an extremal domain of attraction, and $K$ belongs
to the domain of attraction of a probability distribution $G$ (a
notion to be defined precisely), then
$(X,Y)$ follows a CEVM.  

We begin with the \emph{standard case} which means that $(X,Y) \in
[0,\infty)^2$, and $F_Y \in D(G^*_1)$, 
\begin{equation} \label{eqydgstd}
 tF_Y(t\,\cdot) \vconv \nu_1(\cdot) \mt{in} \mplus(0,\infty] \mt{as} t\goesto\infty, 
\end{equation}
where $\nu_1(x,\infty]=x^{-1},\,x>0$ (a formulation equivalent to \eqref{eq:tailconv}) and $K\in D(G)$ meaning
\begin{equation} \label{eqkernconvstd}
 \mckern{t}{t\cdot\,} \wc G(\cdot) \mt{on} [0,\infty].
\end{equation} In what follows, $\xi $ will always be a random
variable with distribution $G$.


\subsection{Standard CEVM Properties}
\label{seccevmstd}

Conditions \eqref{eqydgstd} and \eqref{eqkernconvstd} imply
 $(X,Y)$ follows a CEVM, provided $G \neq \pp{0}$, \ie unit mass at $\{0\}$. 

\begin{thm} \label{thmcevstd}
Suppose that the joint distribution of the random vector $(X,Y)$ on $[0,\infty)^2$ satisfies \eqref{eqydgstd} and \eqref{eqkernconvstd}, where $G$ is a probability distribution on $[0,\infty)$.  Then 
\begin{equation} \label{eqcevstd}
 t\P\big[(X,Y) \in t\cdot\, \big] \vconv \mu(\cdot) \mt{in} \mplus([0,\infty]\times(0,\infty]), 
\end{equation}
 with limit measure $\mu$ given for  $x,y>0,\,\xi\sim G$ by 
\begin{equation} \label{eqlimmeasstd}
 \mu([0,x]\times (y,\infty]) 
=\int_y^\infty G(x/u)\nu_1(du)
=\frac{1}{x}  \int_0^{x/y} G(u) du =
\inv{y}\P[\xi \leq \frac xy ] - \inv{x} \EP
\xi \ind{\xi \leq x/y}.
\end{equation}
Furthermore, $\mu$ satisfies the conditional non-degeneracy conditions \eqref{eqcondnondegennew} provided 
$G \neq \pp{0}$. 
\end{thm}

\begin{pf}
The convergence \eqref{eqcevstd} is special case of 
Proposition 5.1 of \cite{resnick2011asymptotics}; it is an elaboration
of the continuous mapping theorem.
From \eqref{eqlimmeasstd},  $\mu\([0,x]\times(y,\infty]\)$ is continuous in $x$, and not constant provided $G \neq \pp{0}$. 
Also, since $\mu((x,\infty]\times (y,\infty]) = \int_{(y,\infty]}	 \nu_1(du) \P[\xi > x \inv{u} ]$, that $\mu(\{\infty\} \times (y,\infty])=0$ follows from the fact that $G(\{\infty\}) = 0$. 
Therefore, $\mu$ satisfies \eqref{eqcondnondegennew}. 
\end{pf}

\subsubsection{Properties of the limit measure $\mu$.}\label{subsubsec:propsMu}
From \eqref{eqlimmeasstd} we see $\mu$ is continuous in $x$ and $y$ and 
 if $G$ has a density, then so does $\mu$. 
Continuity in \eqref{eqlimmeasstd} holds even if $G$ is degenerate, \ie $G = \pp{c}$ for some $c>0$; see Example \refp{exdiscrete}. Non-degeneracy of $G$  only becomes relevant in the non-standard case. 
Moreover, $\mu$ cannot be a product measure \cite[Lemma
3.1]{das2011conditioning}.
 
From \eqref{eqlimmeasstd} we also observe that the $y$-axis through
the origin is assigned mass proportional to $\Gz$ since  $\mu(\{0\}\times
(y,\infty]) = \inv{y}\Gz$. Mass on vertical slices of space  
depends on $\EP \xi$, since $\mu((x,\infty] \times (0,\infty]) =
\inv{x} \EP\xi \leq \infty$.  
In terms of conditional distributions, \eqref{eqcevstd} implies 
\begin{equation}\label{eq:H}
 \P[ X \leq tx \bgiven Y>t ] \wc H(x) := \mu([0,x]\times (1,\infty])
= \frac{1}{x} \int_0^{x} 
G(u)du.
\end{equation}

\subsubsection{Extending to a larger cone.}
Convergence \eqref{eqcevstd} extends to standard regular variation on
the larger cone $[0,\infty]^2\backslash \{\bz\}$, so that the
distribution of  $(X,Y)$ is in a
bivariate domain of attraction, if and only if  
$F_X \in D(G_1^*)$
as well \cite[Proposition 4.1]{das2011conditioning}. 
In this case, 
\begin{equation} \label{eqmdoastd}
 t\P\big[\inv{t}(X,Y) \in \comp{[\bz,(x,y)]} \big] \conv
 \frac{1}{x}\bigg( 1+ \int_0^{x/y} 
G(u)du \bigg), 
\end{equation}
implying that $\EP \xi \leq 1$, and the $x$-axis receives mass according to $\mu((x,\infty] \times \{0\}) = \inv{x}(1-\EP \xi)$.

\subsubsection{Degenerate $G$; asymptotic independence.}
If $G=\pp{0}$, then the convergence \eqref{eqcevstd}  holds with limit
measure  $\mu([0,x]\times(y,\infty]) = \inv{y}$ but
conditional non-degeneracy \eqref{eqcondnondegennew}
fails, since all the mass lies on the
$y$-axis, so $(X,Y)$ does not follow a standard CEVM.  
This is in fact a manifestation of asymptotic independence. 
Indeed, 
\[ \P[ X>tx \bgiven Y>t ] \goesto 0 \] 
	for any $x$, so, given that $Y$ is extreme (exceeding the threshold $u(t)=t$), it is very unlikely to observe $X$ to be similarly extreme. 
If the joint distribution of $(X,Y)$ is regularly varying on the larger cone $[0,\infty]^2\backslash \{\bz\}$, 
then 
\[ t\P\big[\inv{t}(X,Y) \in \comp{[\bz,(x,y)]} \big] \conv \inv{x} + \inv{y}, \]
which means that $X$ and $Y$ are asymptotically independent in the
usual sense \cite[Section 5]{heffernan2007limit}. In this case,
$(X,Y)$ does not follow a standard CEVM because of degeneracy,
although a CEVM may hold if $X$ is normalized differently; see
Section \ref{sectiongensc}.  

This suggests viewing the parameter $\Gz$ as a measure of  asymptotic dependence from $Y$ to $X$. For example, given $Y$, we could write $X$ as a mixture 
\begin{equation}\label{eqasindep}
X = W X_0 + (1-W) X_1, 
\end{equation}
where $X_0$ and $Y$ are asymptotically independent, $X_1$ and $Y$ are
asymptotically dependent, and $W\sim \text{Bernoulli}(\Gz)$. This is
suggested by  the canonical form of the update function
representation of $K\in D(G)$ \cite[Section 2.3]{resnick2011asymptotics}. Asymptotic
dependence in the reverse direction, given large $X$, would then be
quantified by $1-\EP\xi$ if appropriate. The latter phenomenon is
hinted at by Segers \cite{segers2007multivariate} in his definition of
the ``back-and-forth tail chain'' to approximate stationary Markov
chains .


\subsection{Examples}

Examples illuminate properties of the CEVM based on Markov kernels as
in \eqref{eqkernconvstd}. First, as in \cite[Example
8]{das2011conditioning},
  given any distribution $G$ on
$[0,\infty)$,
 we  construct a CEVM whose limit measure $\mu$ is built on $G$ as
 in  \eqref{eqlimmeasstd}.

\begin{exa} \label{exspecg}
Take $G$ to be any probability distribution on $[0,\infty)$
not concentrating at $0$.
Let $Y\sim\text{Pareto}(1)$ on $[1,\infty)$, $\xi\sim G$, independent of $Y$, and put $X=\xi Y$. 
A version of the conditional distribution is 
\[ \mckern{y}{\cdot} = \P[X \in \cdot \bgiven Y=y ] = \P[ \xi Y
\in\cdot \bgiven Y=y] = G(\inv{y} \cdot), \] 
and  $K$ satisfies
\eqref{eqkernconvstd} and in fact
 and $\mckern{t}{t\cdot } = G(\cdot)$.   
Consequently, $(X,Y)$ follows a standard CEVM with limit measure as in \eqref{eqlimmeasstd}. 
In fact, for $x,y>0$, we have 
\ba
 \P[X\leq x, Y>y] &= \int_{(y,\infty]} \mckern{u}{[0,x]} P[Y\in du] \\ 
 &= \int_{y\smax 1}^\infty \P[\xi \leq x\inv{u}] u^{-2} du =
 \frac{1}{x} \int_0^{x\smin \frac{x}{y}} G(u)du.
\ea
Furthermore, $(X,Y)$ belong to a standard bivariate domain of attraction  \eqref{eqmdoastd} iff $F_X \in D(G_1^*)$ as well.
The marginal distribution of $X=\xi Y$ is 
\[ F_X(x)  =\frac{1}{x}
\int_0^{x} G(u) du = H(x), \]
(from \eqref{eq:H}) which has density $f_X(x) = \inv{x}(G(x) - H(x) \big)$ for $x\geq 0$. Since 
\[ \lim_{t\goesto\infty} t\P[X>tx] = \lim_{t\goesto\infty} \frac{1}{x} \int_0^{tx} \P[\xi > u] du = \inv{x} \EP \xi \ \ (\leq \infty), \]
$(X,Y)$ belongs to the standard domain of attraction iff $\EP \xi=1$.  \qed
\end{exa}

Using the Example \ref{exspecg}
recipe, we  explore the CEVM in a variety of special cases. 

\begin{exa}
Choose $\xi\sim \text{Exp}(\lambda)$ and we have 
$X = \inv{\lambda}{Y}E$, where $E\sim\text{Exp}(1)$. 
The limit measure is 
\[ \mu([0,x]\times(y,\infty]) 
= \frac{1}{x} \int_0^{x/y} (1-e^{-\lambda u}) du = \frac{1}{y}
-\frac{1}{\lambda x} + \frac{e^{-\lambda x/y}}{\lambda x}, \] 
and 
the marginal distribution of $X$ is  $F_X(x) = 1 - \inv{(\lambda x)}(1-{e^{-\lambda x}})$
with density $f(x) = \inv{\lambda} x^{-2} (1-e^{-\lambda x}) -
\inv{x}e^{-\lambda x}$, and $F_X$ satisfies \eqref{eqydgstd} iff
$\lambda=1$.  \qed
\end{exa}

Next, we suppose $\xi$ is heavy-tailed. 

\begin{exa}
For $\alpha>0$ let $\xi\sim \text{Pareto}(\alpha)$ so $1-G(x)=:\bar
G(x)=x^{-\alpha},\, x\geq 1$. 
The limit measure assigns no mass to $\{(x,y) : 0 \leq x \leq y\}$, and for $x>y>0$, 
\[ 
\mu([0,x]\times(y,\infty])
 = \begin{cases}
\displaystyle \frac{1}{y} - \(\frac{\alpha}{\alpha-1}\)\frac{1}{x} + \frac{y^{\alpha-1}}{x^{\alpha}(\alpha-1)} & \alpha > 1\vspace{2mm} \\
\displaystyle \frac{1}{y} - \frac{1}{x} - \frac{\log x}{x} + \frac{\log y}{x} & \alpha = 1 \vspace{2mm} \\
\displaystyle \frac{1}{y} + \(\frac{2-\alpha}{1-\alpha}\)\frac{1}{x} + \frac{1}{x^{\alpha}y^{1-\alpha}(1-\alpha)} & \alpha < 1 .
\end{cases}
\]
When $\alpha \leq 1$, $\EP\xi = \infty$  and 
$\mu((x,\infty]\times(y,\infty]) = {\inv{y} -
  \mu([0,x]\times(y,\infty])} \goesto\infty$  
as $y\downarrow 0$. \qed
\end{exa}

It is also possible that $G$ is discrete, although the CEVM limit measure $\mu$ remains continuous. 

\begin{exa} \label{exdiscrete}
Suppose $\xi$ has discrete distribution $\P[\xi=k] = a_k$, $k=0,1,\dots$. 
In this case, the limit measure is given by 
\[ \mu([0,x]\times(y,\infty]) = \frac{1}{x} \int_0^{x/y} \Big(\sum_{k=0}^{[u]} a_k \Big) du 
= \sum_{k=0}^{[x/y]} a_k (\inv{y}- k \inv{x}), 
\]
which is continuous in $x$ and $y$, and $F_X(x) = \sum_{k=0}^{[x]} a_k (1- k \inv{x})$. 
In particular, if $\P[\xi=c] = 1$ for some $c>0$, we obtain 
\[ \mu([0,x]\times(y,\infty]) = (\inv{y}- c \inv{x})\ind{x>cy>0}. \] 
The conditional non-degeneracy conditions \eqref{eqcondnondegennew}
are satisfied even though $G$ is degenerate.  \qed
\end{exa}

The final example shows how $G$ reflects 
asymptotic independence between $X$ and $Y$.

\begin{exa}\label{exdegen}
Consider $Y\sim$ Pareto$(1)$, and $Z$  
independent of $Y$ such that $\P[Z<\infty] = 1$. 
Take $ X=Y\smax Z .$
Given $Y$ is extreme, it is unlikely that $Z$ is as extreme as $Y$ since they are independent. 
We have 
\[ \mckern{y}{[0,x]} = \P[Y \leq x, Z\leq x \given Y=y] = \P[Z \leq x] \ind{x \geq y}, \]
and so 
\[ \mckern{t}{t[0,x]} = \P[Z \leq tx]\ind{x \geq 1} \conv 
\ind{x\geq 1} =
\pp{1}([0,x])  = G([0,x]). \]
As in the previous example, the limit measure is 
\[ \mu([0,x]\times(y,\infty]) = (\inv{y}-  \inv{x})\ind{x>y>0}. \] 
On the other hand, if $X'=Y\smin Z $ then when $Y$ is large,  it is likely $X' = Z$, so $X'$
should be asymptotically independent of $Y$. More precisely,
\[ \mckern{y}{(x,\infty]} = \P[Y>x, Z>x \given Y=y] = \P[Z>x] \ind{y>x}, \]
from which 
\[ \mckern{t}{t(x,\infty]} = \P[Z>tx] \ind{x<1} \conv 0 \] 
for $x>0$. 
Therefore, $G = \pp{0}$, and the conditional non-degeneracy conditions
fail.  \qed
\end{exa}

\subsection{Counter-examples}

As expected, the converse to Theorem \ref{thmcevstd} can fail.
If $(X,Y)$ follows a non-degenerate CEVM as in \eqref{eqcevstd}, and $K$ is a specific version of the conditional distribution $\P[X\in \cdot \given Y=y]$, it does not necessarily follow that there exists a distribution $G$ such that \eqref{eqkernconvstd} holds. 
The failure of  \eqref{eqkernconvstd} can happen  in two ways. 
There may  exist a probability distribution $G$ on $[0,\infty]$ satisfying \eqref{eqkernconvstd} 
with ${G(\{\infty\})>0}$ or
it may be possible to obtain two distinct limit distributions down different subsequences 
$\{t_n\}$ and $\{t'_n\}$. 

Example \ref{exdefecg}
 where $G(\{\infty\})>0$
 emphasizes the importance of assuming $\mu(\{\infty\} \times (y,\infty])=0$. 

\begin{exa} \label{exdefecg}
As usual, take $Y\sim$ Pareto$(1)$ and suppose that 
\[ X = WY + (1-W) Y^2, \]
 where $W\sim\text{Bernoulli}(p)$ independent of $Y$.  
Then 
\[\mckern{y}{\cdot} = \P[X\in \cdot \given Y=y] = p \pp{y} + (1-p) \pp{y^2}, \] 
so 
\[ \mckern{t}{t\,\cdot} = p \pp{1} + (1-p) \pp{t} \wc p \pp{1} + (1-p) \pp{\infty} = G \mt{on} [0,\infty]. \]
Indeed, 
for $0 \leq x < \infty$, 
\[ \mckern{t}{t[0,x]}  = p\pp{1}([0,x]) + (1-p) \pp{t}([0,x])
\conv p \pp{1}([0,x]) 
\]
showing that $G(\{\infty\}) = 1-p$.

On the other hand, for $x,y>0$, 
\ba
\P[X \leq x,\, &Y>y] = p\P[Y\leq x,\, Y>y] + (1-p)\P[Y^2\leq x,\, Y>y] \\
 &= p\bigg[\frac{1}{(y\smax 1)}-\frac{1}{x}\bigg] \ind{x \geq (y\smax 1)} + (1-p) \bigg[\frac{1}{(y\smax 1)} - \frac{1}{\sqrt{x}} \bigg] \ind{x \geq (y\smax 1)^2}, 
\ea 
so for $t$ sufficiently large, 
\begin{align}
t\P[X\leq t x, Y> t y] &= p\bigg[\frac{1}{y}-\frac{1}{x}\bigg] \ind{x \geq y} + (1-p) \bigg[\frac{1}{y} - \frac{\sqrt{t}}{\sqrt{x}} \bigg] \ind{\sqrt{x}/\sqrt{t} \geq y} \nonumber  \\
 &\conv  p\bigg[\frac{1}{y}-\frac{1}{x}\bigg] \ind{x \geq y} = \mu([0,x] \times (y,\infty]).  \label{eqexdefec}
\end{align}
The measure $\mu$ assigns positive mass to  $\{\infty\} \times
(y,\infty]$ since
\ba 
\mu((x,\infty] \times (y,\infty]) = \inv{y} - \mu([0,x] \times (y,\infty]) = \frac{1}{y}\ind{x<y} + \bigg[\frac{1-p}{y} + \frac{p}{x} \bigg] \ind{x \geq y}, 
\ea 
and thus $\mu(\{\infty\} \times (y,\infty]) = (1-p)\inv{y}$. Therefore, $\mu$ does not satisfy \eqref{eqcondnondegennew}. 

Under a different normalization, we  obtain a proper limit $G$. 
Indeed, note that 
\[ \mckern{t}{t^2 \cdot}  = p\pp{\inv{t}} + (1-p) \pp{1} \wc p \pp{0} + (1-p) \pp{1} \sim \text{Bernoulli}(1-p), \]
and hence, 
\ba 
t\P[X\leq t^2 x, Y> t \cdot y] = &p(\inv{y}-\inv{(tx)}) \ind{x \geq y/t} + (1-p) (\inv{y} - x^{-1/2}) \ind{x \geq y} \\ 
 &\conv  p\inv{y}+(1-p)(\inv{y}-x^{-1/2})\ind{x \geq y}. 
\ea
This limit does satisfy \eqref{eqcondnondegennew}. \qed
\end{exa}

 Without
condition  \eqref{eqcondnondegennew}, the convergence of types theorem
fails and it is possible to obtain different CEV limits under different normalizations. 
From  \eqref{eqlimmeasstd},  $\mu(\{\infty\} \times (y,\infty]) =
G(\{\infty\}) \inv{y}$
and  excluding defective distributions in Theorem \ref{thmcevstd}  avoids cases like the previous one.

Here is an example of a CEVM 
 where the normalized kernel $K$ does not have a unique limit.

\begin{exa} \label{exnotconv}
Suppose $Y\sim$ Pareto$(1)$, and define $X$ by 
\[ X = W Y + (1-W) 2Y\ind{Y \in 
[0,\infty) \backslash\N}
\] 
where $W\sim\text{Bernoulli}(p)$ independent of $Y$. 
In other words, given $Y=y$, $X$ takes the value $y$ or $2y$ according
to a coin flip, unless $y$ is an integer, in which case $X$ will be
either $y$ or 0.  The CEVM holds for  $(X,Y)$. Since $\P[Y \in \N] = 0$, we have  
\ba 
\P[X\leq x, Y>y] &= \P[X\leq x,\, Y>y,\, Y\in [0,\infty) \backslash\N ] \\
 &= p\P(Y \leq x, Y>y) + (1-p) \P(2Y\leq x, Y>y) \\
 &= p(\inv{y}-\inv{x})\ind{x\geq y} + (1-p) (\inv{y}-2\inv{x})\ind{x\geq 2y}, 
\ea 
and $t\P[X\leq tx, Y>ty] = \P[X\leq x, Y>y]$, which satisfies
\eqref{eqcevdisth}
and the requirement that $\mu( (\cdot) \times (y,\infty])$ not be
degenerate for any $y$.
However, the conditional distribution of $X$ given $Y$ is 
\[ \mckern{y}{\cdot} = 
\begin{cases} 
 p \pp{y} + (1-p) \pp{0} & y \in \N   \\ 
 p \pp{y} + (1-p) \pp{2y} & y \in [0,\infty) \backslash \N
\end{cases}, \]
so 
\[ \mckern{t}{t\,\cdot} = 
\begin{cases} 
p \pp{1} + (1-p) \pp{0} & t \in \N   \\ 
 p \pp{1} + (1-p) \pp{2} & t \in [0,\infty) \backslash \N
\end{cases}. \]
We obtain different limits along the sequences $t_n = n$ and $t'_n =
n/2$
and $\mckernl{t}{t\,\cdot}$ does not converge. \qed
\end{exa}

The technical difficulty highlighted in 
Example \ref{exnotconv}  is that conditional distributions of the form $\P[X\in
\cdot \given Y=y]$ are only specified up to sets of $\P[Y\in
\cdot\,]$-measure zero. If $Y$ is absolutely continuous, we can
alter the conditional probability for a countable number of $y$
without affecting the joint distribution.  Consequently, constructing
 a convergence theory based on conditional
distributions requires care. The best one can do is fix a version of the kernel or,
if circumstances allow, choose a version of the kernel with some claim
to naturalness based on smoothness.
This is the reason the approach in \cite{
das2011conditioning,
heffernan2007limit} is based on vague convergence of measures rather
than convergence of conditional distributions as in \cite{heffernan2004conditional}.


\section{General Normalization for $X$}
\label{sectiongensc}
The CEVM allows different  normalizations for $X$ and $Y$, as in
\eqref{eqcevmdef}, but the formulation  $K\in D(G)$  in
\eqref{eqkernconvstd}, imposes the same normalization for both.  
We now allow general linear normalizations of $X$ in the kernel
condition, continuing to assume condition \eqref{eqydgstd} 
that $F_Y \in D(G^*_1)$.

We will assume the following generalization of \eqref{eqkernconvstd}:
there exist scaling and centering 
functions $\alpha (t)>0$, $\beta (t) \in
\R$ and 
 a non-degenerate probability distribution $G$ on
$[-\infty,\infty)$,  such that  
\begin{equation} \label{eqkernconvgen}
\mckernb{t}{[-\infty,\alpha(t) x + \beta(t)]} \wc G([-\infty,x]) \mt{on} [-\infty,\infty].  
\end{equation}


\subsection{CEVM Properties}

\label{seccevmgen}

Consider the decomposition 
\begin{equation}\label{eq:probConv}
t\P\bigg[ \frac{X-\beta(t)}{\alpha(t)}\leq x, Y > ty \bigg] = \int_{(y,\infty]} t\P[Y \in t du ] \mckernb{tu}{ [-\infty,\alpha(t)x + \beta(t)]}. 
\end{equation}
By a variant of the continuous mapping theorem (\cite[Lemma
8.2]{resnick2011asymptotics}, the integrals  converge provided  
$\mckernb{tu(t)}{ [-\infty,\alpha(t)x + \beta(t)]} \goesto
\varphi_x(u)$ whenever $u(t) \goesto u>0$.  Proposition
\ref{propkernunifgen} 
discusses when  this happens.

\newcommand{\tkg}[1][G]{\kappa_{#1}}

Given $\rho,k \in \R$, 
define the \emph{generalized tail kernel} associated with a distribution $G$ on $[-\infty,\infty]$ as the transition function $\kappa_G : (0,\infty) \times \Bor[-\infty,\infty] \goesto[0,1]$ given by 
\begin{equation} \label{eqtailkerngen}
 \mckern[\tkg]{y}{A} = 
G\big(y^{-\rho}[A-\psi(y)] \big), 
\end{equation}
where $\psi$ is specified in \eqrefp{eqervpsi}. 
Note that $\tkg$ describes transitions between two different spaces. 
Since $\psi$ satisfies $\psi(uy) = u^{\rho} \psi(y) + \psi(u)$,  a
kernel $\kappa$ has the form \eqref{eqtailkerngen} iff 
\begin{equation}\label{eq:scaling}
\mckern[\kappa]{uy}{A} = \mckernb[\kappa]{y}{u^{-\rho}[A-\psi(u)]}. \end{equation}

\begin{prop} \label{propkernunifgen}
Let $K : (0,\infty) \times \Bor[-\infty,\infty] \goesto[0,1]$ be a transition function satisfying \eqref{eqkernconvgen} with $G$ is non-degenerate.  
There exists a family of non-degenerate probability 
distributions $\{G_u : 0 < u < \infty\}$ on $[-\infty,\infty)$ such
that  for $0 < u < \infty, $
\begin{equation} \label{eqpropkernconvfam}
\mckernb{tu}{[-\infty,\alpha(t) x + \beta(t)]} \wc
G_u([-\infty,x])\quad \text{on } [-\infty,\infty], 
\end{equation} 
as $t\goesto\infty$ if and only if $\alpha, \beta \in \ERV_{\rho,k}$ 
as in \eqrefp{eqdeferv}. 
In this case, $G_1 = G$, and 
\begin{equation} \label{eqkernunifgen}
\mckernb{tu_t}{[-\infty,\alpha(t) x + \beta(t)]} \wc 
\mckern[{\tkg[{G}]}]{u}{[-\infty,x]} \mt{on} [-\infty,\infty]
\end{equation} 
whenever $u_t = u(t) \goesto u \in (0,\infty)$; i.e., the limit is a transition function of the form \eqref{eqtailkerngen}, where 
$\rho,k$ are the ERV parameters of $\alpha,\beta$. 
\end{prop}

\begin{pf}
Assume first that $\alpha, \beta \in \ERV_{\rho,k}$ and define
\[ h_t(y;u) = \frac{\alpha(tu)}{\alpha(t)} y + \frac{\beta(tu)-\beta(t)}{\alpha(t)},\]
 so that by \eqref{eqdeferv}, $h_t(y_t;u) \goesto h(y;u) = u^{\rho}y+\psi(u)$ whenever $y_t \goesto y \in \R$. 
For $u>0$, 
\[ \mckernb{tu}{[-\infty,\alpha(t) x + \beta(t)]} = 
\mckernb{tu}{\alpha(tu)\{\inv{h_t}(\cdot\,;u)[-\infty, x]\} + \beta(tu)}.
\] 
Applying the second continuous mapping theorem
(\cite{billingsley1968convergence}, \cite[Lemma 8.1]{resnick2011asymptotics})
to 
\eqref{eqkernconvgen}, we have 
\[ \mckernb{tu}{[-\infty,\alpha(t) x + \beta(t)]} \wc (G \circ
\inv{h}(\cdot\,;u)) ([-\infty, x])=
G([-\infty,(x-\psi(u) )/u^\rho]).\] 
Hence, \eqref{eqpropkernconvfam} holds with $G_u (\cdot)=
\mckernl[\tkg]{u}{\cdot}$
 and  $G_1 = G$. 
Furthermore, we have $h_t(x_t ; u_t) \goesto h(x ; u)$ whenever ${u_t\goesto u>0}$, establishing \eqref{eqkernunifgen}.

For the converse, we employ convergence of types. Denote by 
$H_t(\cdot)$ the distribution $\mckernl{t}{\cdot}$. 
Then, on the one hand, we have $H_t([-\infty,\alpha(t) x + \beta(t)]) \wc G_1([-\infty,x])$. On the other hand, fixing $c>0$, we have 
\[
H_t(\alpha(tc) x + \beta(tc)) = \mckernb{(tc)\inv{c}}{[-\infty,\,\alpha(tc) x + \beta(tc)]} \wc G_{\inv{c}}([-\infty,x]).
\] 
Convergence of types yields that $\alpha,\beta \in \ERV_{\rho,k}$, and 
\[ G_{\inv{c}}([-\infty,x]) = G_1([-\infty,c^{\rho} x + \psi(c)]), \] with $\psi$ as in \eqref{eqervpsi}. 
Using the identity \eqref{eqervpsiinv}, 
we find that $G_u$ has the form \eqref{eqtailkerngen}, with $G=G_1$.
\end{pf}

Starting from kernel convergence
\eqref{eqkernconvgen},
Proposition \ref{propkernunifgen} implies that 
$\alpha, \beta$ being ERV is
 necessary and sufficient for obtaining a CEVM.
Unlike Section \ref{stdcontbasiccase}, here we need
 $G$ to be non-degenerate in order to apply
the convergence of types theorem.

\begin{thm} \label{thmcevmgen}
Suppose $(X,Y)$ is a random vector on $\R\times [0,\infty)$ and
\eqref{eqydgstd} holds.
Assume  \eqref{eqkernconvgen} holds for 
non-degenerate limit distribution $G$
on $[-\infty,\infty)$ and  scaling and centering functions
$\alpha (t) >0$ and $\beta (t)\in\R$.
As $t\goesto\infty$,  
\begin{equation} \label{eqcevgen}
 t\P\bigg[\bigg(\frac{X-\beta(t)}{\alpha(t)},
 \frac{Y}{t}\bigg)\in\cdot\,\bigg] \vconv \mu(\cdot) \not \equiv 0 \mt{in}
 \mplus([-\infty,\infty]\times(0,\infty])  
\end{equation}
where $\mu$  satisfies \eqref{eqcondnondegennew},
if and only if $\alpha, \beta \in \ERV_{\rho,k}$.  
In this case,  $\mu$ is specified by 
\begin{equation} \label{eqproplimmeasgen}
\mu\([-\infty,x]\times(y,\infty]\) = \int_{(y,\infty]} \nu_1(du)
G\bigl( \rinv{u}(x-\psi(u)) \bigr) , \quad x\in\R,\ y>0,  
\end{equation}
with $\psi$ as in \eqref{eqervpsi} and $\nu_1(du)=u^{-2}du,\,u>0$. Expression
\eqref{eqproplimmeasgen} is continuous in $x$ and $y$ if $(\rho,k)
\neq (0,0)$, or if $G$ is continuous.  
\end{thm}

\begin{pf}
The convergence \eqref{eqcevgen} to a limit $\mu$ satisfying \eqref{eqcondnondegennew} implies $\alpha,\beta\in\ERV$ \cite[Proposition 1]{heffernan2007limit}. 
Conversely, if $\alpha, \beta \in \ERV_{\rho,k}$, then the convergence
\eqref{eqcevgen} follows from 
a variant of the continuous mapping theorem  (\cite[Lemma
8.4]{resnick2011asymptotics})
in light of 
\eqref{eqkernunifgen}, yielding the limit in \eqref{eqproplimmeasgen}.  
We check that $\mu([-\infty,x]\times(y,\infty])$ is continuous when $(\rho,k)\neq (0,0)$
by applying dominated convergence: if $x_n\goesto x$, then 
\[ G\bigl( \rinv{u}(x_n-\psi(u))\bigr) \goesto G\bigl(\rinv{u}(x-\psi(u))\bigr)\] 
for all except a countable number of $u$ corresponding to discontinuities of the distribution function. Continuity in $y$ is clear. 
Also, if $(\rho,k)= (0,0)$, then $\mu([-\infty,x]\times(y,\infty]) =
\inv{y}G(x)$, which is continuous if $G$ is. In either case,
$\mu([-\infty,x]\times(y,\infty])$ is non-degenerate in $x$ because
$G$ is non-degenerate.  
Finally, $\mu(\{\infty\} \times (y,\infty]) = \inv{y}G(\{\infty\}) = 0$. Therefore, $\mu$ satisfies \eqref{eqcondnondegennew}. 
\end{pf}

Changing variables $u\mapsto 1/u$ in \eqref{eqproplimmeasgen},
the limit measure is 
\begin{equation} \label{eqlimmeasgen}
\mu([-\infty,x]\times(y,\infty]) = \int_0^{\inv{y}} G( u^\rho x+\psi(u) )\,du , 
\end{equation}
where 
\[ u^\rho x+\psi(u) = \begin{cases} u^\rho(x+k\inv{\rho}) - k\inv{\rho} & \rho \neq 0 \\ x + k\log u & \rho = 0 \end{cases}. \]
Changing variables, we obtain the following expressions for $\mu$ according to $(\rho,k)$: 
\begin{align} \label{eqlimmeasgencases}
 &\mu([-\infty,x]\times(y,\infty]) = \\
&\begin{cases} 
\displaystyle\frac{1}{\rho \abs{x+k\inv{\rho}}^{1/\rho}}
\int_{0}^{\abs{x+k\inv{\rho}}\rinv{y}} u^{(1-\rho)/\rho}  G( u \sgn
(x+k\inv{\rho}) - k\inv{\rho}) du & \rho \neq 0 \vspace{3mm} \\  
 \displaystyle\frac{1}{\abs{k} e^{x/k}}\int_{-\infty}^{x\sgn(k)
   -\abs{k}\log y} e^{u/\abs{k}} G( u\sgn(k)) du
 \hspace{\stretch{1}} \rho = 0, & 
k \neq 0 \vspace{3mm} \\
\inv{y}G( x)	\hspace{\stretch{1}} \rho = 0, & 
k=0
\end{cases}. \nonumber
\end{align}
Here $\sgn(v) = v/\abs{v} \ind{v \neq 0}$, and we 
read the measure as $\inv{y}G( -k\inv{\rho})$ when $x=-k\inv{\rho}$ for the case $\rho \neq 0$. 
Continuity in $x$ and $y$ when $(\rho,k) \neq (0,0)$ is apparent from the above expressions. 

We  give an example where $K$ satisfies \eqref{eqkernconvgen}, but \eqref{eqcevgen} fails because $\alpha, \beta$ are not ERV. 

\begin{exa} \label{excevmgen}
Consider $Y\sim \text{Pareto}(1)$ and $U \sim \text{Uniform}(0,1)$, independent of $Y$. Put $X = Ue^Y$. Then 
\[ \mckern{y}{[0,x]} = \P[X\leq x \given Y=y] = \P[U\leq xe^{-y}] = xe^{-y} \smin 1. \]
Polynomial scaling is not strong enough to give an informative limit,
since  for any $\rho>0$,
\begin{equation*} 
\mckern{t}{t^\rho [0,x]} = x^\rho t^\rho e^{-t} \smin 1 \goesto 0. 
\end{equation*}
The appropriate normalization is exponential $(\alpha(t),\beta(t)) =(
e^t,0)$, which is not ERV:
\[ \mckern{t}{\alpha(t) [0,x]} = x e^t e^{-t} \smin 1 \goesto x \smin 1 = G(x), \]  
and  up to asymptotic
equivalence,  by the convergence to types theorem, this the only
normalization yielding a non-degenerate limit.
Since $(\alpha(t), \beta(t))$ is not ERV, Theorem \ref{thmcevmgen}
claims that $(X,Y)$ cannot follow a CEVM. To verify this {\it ab
  initio\/}, consider $y>0$ and $t$ so large that $ty>1$:
\bml
 t\P[X \leq \alpha(t) x , Y>ty ] =tP[Ue^Y\leq e^tx,Y>ty]=
 \int_y^\infty u^{-2} \bigl(xe^{-t(u-1)}\wedge 1\bigr)    du                   \\
=  \int_{((y\smax 1),\infty)} u^{-2}du  \{ xe^{-t(u-1)}  \smin 1 \} 
 + \ind{y<1} \int_{(y,1]} u^{-2}du  \{ xe^{-t(u-1)}  \smin 1 \} . 
\eml
The first integral in the previous sum is bounded by $x \inv{y}  e^{-t(y-1)} \goesto 0$. If $y\leq 1$, the second integral approaches $\nu_1(y,1] = \inv{y}-1$. Therefore, the limit is degenerate in $x$, violating conditional non-degeneracy \eqref{eqcondnondegennew}. 
This is not repaired by using ERV normalization since if
 $\tilde{\alpha}, \tilde{\beta}$ are ERV,  then 
\ba
 t\P[X \leq \tilde{\alpha}(t) x + \tilde{\beta}(t),\, &Y>ty ] =
 \int_{(y,\infty )} u^{-2} du  \mckernb{tu}{[0,\tilde{\alpha}(t) x + \tilde{\beta}(t)]} \\ 
 &= \int_{(y,\infty )} u^{-2}du \big\{ e^{-tu}(\tilde{\alpha}(t) x + \tilde{\beta}(t)) \smin 1 \big\} \conv 0, 
\ea
which follows from the asymptotic properties of ERV functions (see
Section \refp{secerv}). \qed \end{exa} 


\subsection{Standardization of $X$}

In certain cases, it is possible to standardize the $X$ variable
\cite[Section 3.2]{das2011conditioning}.
\subsubsection{Standardization functions.}\label{subsubsec:stanFcts}
Denote by $x^*$ and $x_*$ the upper and lower endpoints of the distribution of $X$ respectively, \ie 
\[ x^* = \sup \{ x : F_X(x)<1 \} \mand x_* = \inf \{ x : F_X(x)>0 \}. \]
Call $f:(0,\infty) \mapsto (x_*,x^*)$ a \emph{standardization
  function} if $f$ is monotone and $\lim_{x\goesto\infty} f(x) =x^*$
if $f$ is non-decreasing and $\lim_{x\goesto\infty} f(x) =x_*$ if $f$ is non-increasing.
As in  \cite[Section 3]{das2011conditioning}, we standardize with such
functions.
For the purpose of this section, extend the definition of
$\ginv{f}:(x_*,x^*)\mapsto (0,\infty)$ in order to invert right-continuous monotone functions which are either increasing or decreasing. 
Define 
\[ \ginv{f}(x) = \begin{cases} \inf\{ y : f(y) \geq x \} & \quad\text{if } f \text{ is non-decreasing} \\ \inf\{ y : f(y) \leq x \} & \quad\text{if } f \text{ is non-increasing} \end{cases}. 
\]
Note that $\ginv{f}$ is left-continuous for $f$ non-decreasing and right-continuous for $f$ non-increasing. 
The main property we shall be using is that 
\[ \left\{
\begin{aligned} 
&\ginv{f}(x) \leq y \iff x \leq f(y) &\qquad& f \text{ non-decreasing} \\ 
&\ginv{f}(x) \leq y \iff x \geq f(y) && f \text{ non-increasing}
\end{aligned}\right. .
\]
The distinction between the two cases is a technicality which should not cause confusion in the following discussion. 
Also,  say that a monotone function $f$ has two {\it points of change\/}
if there exist $x_1 < x_2 < x_3$ such that $f(x_1) < f(x_2) < f(x_3)$
for $f$ non-decreasing, and with the opposite inequalities in the
non-increasing case.  

If the pair $(X,Y)$ satisfies \eqref{eqcevgen} for some $\alpha>0$ and $\beta$,
then we say $(X,Y)$ can be \emph{standardized} if there exists a
standardization function $f$ and a  non-null Radon measure  $\mu^*$  such that  
\begin{equation} \label{eqcevstdx}
 t \P\big[ \inv{t} (\ginv{f}(X),Y) \in \cdot\, \big] \vconv
 \mu^*(\cdot) \mt{in} \mplus([0,\infty]\times (0,\infty]).
\end{equation}
If the limit $\mu$ in \eqref{eqcevgen} satisfies the conditional non-degeneracy conditions \eqref{eqcondnondegennew}, then standardization is possible if and only if $(\rho,k) \neq (0,0)$; \ie $\mu$ is not a product. 

\subsubsection{Characterizing standardization functions.}\label{subsubsec:stan}
What property does $f$ need to be a standardization function
satisfying \eqref{eqcevstdx}?

\nc{\Ax}{A_\varphi(x)}

\begin{prop} \label{propstdfn}
Suppose $(X,Y)$ follow a CEVM so that  \eqref{eqcevgen} holds with $\mu$
satisfying the conditional non-degeneracy conditions
\eqref{eqcondnondegennew}. Assume $(\rho,k) \neq (0,0)$.  
A function $f$ standardizes $(X,Y)$ in the sense of
\eqref{eqcevstdx} where $\mu^*$ satisfies the conditional
non-degeneracy conditions iff 
\begin{equation} \label{eqstdfn}
\frac{f(tx) - \beta(t)}{\alpha(t)} \conv \varphi(x), \qquad x>0, 
\end{equation}
where $\varphi$ has at least two points of change. In this case, $\mu$ and $\mu^*$ are related by  
\[\mu^*([0,x] \times (y,\infty]) = \mu(\Ax \times (y,\infty]), \]
where 
\begin{equation} \label{eqax}
 \Ax = \begin{cases} [-\infty,\varphi(x)] & \quad f \text{ non-decreasing} \\ 
[\varphi(x),\infty] & \quad f \text{ non-increasing} \end{cases}\ . 
\end{equation}
\end{prop}

It follows that $\alpha (\cdot), f(\cdot) \in \ERV$, though not
necessarily with the same parameters as $\alpha,\beta$. However,
depending on the case, $f$ can be expressed in terms of either
${\beta}$ or ${\alpha}$ (\cite[Proposition
2.3.3]{das2009conditional}).  

\begin{pf}
Suppose $f$ is non-decreasing. Then for $x,y>0$, we can write 
\begin{equation} \label{eqstdprob}
t\P\bigg[ \frac{\ginv{f}(X)}{t} \leq x, \frac{Y}{t} > y \bigg] = t\P\bigg[ \frac{X-\beta(t)}{\alpha(t)} \leq \frac{f(tx) - \beta(t)}{\alpha(t)}, \frac{Y}{t} > y \bigg]. 
\end{equation} 
If $f$ satisfies \eqref{eqstdfn}, then \eqref{eqcevstdx} holds with 
\[ \mu^*([0,x] \times (y,\infty]) = \mu([-\infty,\varphi(x)]\times (y,\infty])
\] 
non-degenerate in $x$. 
On the other hand, if \eqref{eqcevstdx} holds, then \eqref{eqstdprob} implies \eqref{eqstdfn}, and $\varphi$ has at least two points of increase because $\mu^*$ is non-degenerate in $x$. The mass at $\{\infty\}$ condition in \eqref{eqcondnondegennew} follows from the fact that $\lim_{x\goesto\infty} \varphi(x) = \infty$ if $f$ is non-decreasing (see \eqref{eqstdfnphi} below). The case for $f$ non-increasing is similar, after reversing the inequality for $X$ on the right-hand side of \eqref{eqstdprob}. 
\end{pf}

Assuming \eqref{eqstdfn} and $\alpha,\beta \in \ERV_{k,\rho}$, write 
\[ \frac{f(tx) - \beta(t)}{\alpha(t)} = \frac{\alpha(tx)}{\alpha(t)}\frac{f(tx) - \beta(tx)}{\alpha(tx)} + \frac{\beta(tx) - \beta(t)}{\alpha(t)} \] 
 and with  $c = \varphi(1)$,  $\varphi$ has the form 
\begin{equation} \label{eqstdfnphi}
\varphi(x) = \begin{cases} c x^\rho + k\inv{\rho}(x^\rho - 1) & \rho \neq 0 \\ c + k \log x & \rho=0 \end{cases}. 
\end{equation}
If $\varphi$ has two points of change, we get the  constraint that $c\neq 0$ if $\rho \neq 0$, $k = 0$.

\subsubsection{Kernel convergence and standardization.}\label{subsubsec:kernel}
Assuming a standardization function exists and that the conditional
distribution of $X$ given $Y$  
satisfies kernel convergence assumption \eqref{eqkernconvgen},
we can  standardize directly through the Markov kernel. We consider
the new kernel $K_f(y,A)=K(y, f(A))=:\P[f^\leftarrow (X) \in A|Y=y]$.
The next result may be compared to the formulation in \cite[Proposition
2.3.3]{das2009conditional} for joint
distributions rather than Markov kernels.

\nc{\Axi}{A_{\ginv{\varphi}}(x)}

\begin{prop} \label{propkernstdx}
Suppose the transition function $K : (0,\infty) \times \Bor[-\infty,\infty] \goesto[0,1]$ satisfies \eqref{eqkernconvgen} 
for a probability distribution $G$ on $[-\infty,\infty)$. 
If $f$ is a monotone function satisfying \eqref{eqstdfn},
then the transition function $K_f : (0,\infty) \times \Bor[0,\infty] \goesto[0,1]$ defined as 
\[ \mckern[K_f]{y}{A} = \mckern{y}{{f}(A)} \]
satisfies  \eqref{eqkernconvstd},
\[\mckern[K_f]{t}{t[0,x]} 
\wc G(\Ax) 
=: G_f([0,x]) \mt{on} [0,\infty], \]
with $\Ax$ as in \eqref{eqax}. 
Conversely, if we start with a kernel $K(y,\cdot)$ satisfying 
\eqref{eqkernconvstd}, for limit probability measure $G$ on $[0,\infty)$,
then given ERV functions $\alpha>0$, $\beta \in \R$ on $(0,\infty)$,
if $f$ is  monotone  on $(0,\infty)$ satisfying
\eqref{eqstdfn},  
the transition function $\overline{K}_f : (0,\infty) \times
\Bor[-\infty,\infty] \mapsto
[0,1]$ given by
\[ \mckern[\overline{K}_f]{y}{A} = \mckern{y}{\ginv{f}(A)} \]
satisfies \eqref{eqkernconvgen},
\[ \mckernb[\overline{K}_f]{t}{[-\infty,\alpha(t) x + \beta(t)]} \wc G(\Axi) 
=: 
\overline{G}_f([-\infty,x]) \mt{on} f([0,\infty]), \] 
where 
\begin{equation*} 
 \Axi = \begin{cases} [0,\ginv{\varphi}(x)] & \quad f \text{ non-decreasing} \\ 
[\ginv{\varphi}(x),\infty] & \quad f \text{ non-increasing} \end{cases}\ . 
\end{equation*}
\end{prop}

\begin{pf}
Assume \eqref{eqkernconvgen} and $f$ is a non-decreasing
 function satisfying \eqref{eqstdfn}.  
Then, 
\begin{align*} 
\mckernb[K_f]{t}{&t[0,x]} = \mckernb{t}{[-\infty,{f}(tx)]}  \\ 
 &= K \bigg(t \,,\, \alpha(t)\bigg[-\infty,\frac{{f}(tx) - \beta(t)}{\alpha(t)}\bigg] + \beta(t) \bigg) 
\wc G([-\infty,{\varphi}(x)]). 
\end{align*}
Conversely, if $f$ satisfies \eqref{eqstdfn} for $\alpha,\beta\in
\ERV$, then inverting \eqref{eqstdfn}
 yields 
\[ {\ginv{f}(\alpha(t)x+\beta(t))}/{t} \conv \ginv{\varphi}(x), \qquad x \in f((0,\infty)). \]
Consequently, 
\[
\mckernb[\overline{K}_f]{t}{[-\infty,\alpha(t) x + \beta(t)]} = \mckernb{t}{t[0,\inv{t}\ginv{f}(\alpha(t) x + \beta(t))]}  
\wc G([0,\ginv{\varphi}(x)]). 
\]
The case for non-increasing $f$ is similar. 
\end{pf}

If $K$ is a version of the conditional distribution ${\P[X\in \cdot
  \given Y=y]}$ satisfying \eqref{eqkernconvgen}, where $\alpha,\beta 
\in \ERV_{\rho,k}$ with $(\rho,k) \neq (0,0)$ and $F_Y$ is in the
standardized domain of attraction,  
then $(X,Y)$ follows a CEVM by Theorem \ref{thmcevmgen}. 
Furthermore, $(X,Y)$ can be standardized in the sense of
\eqref{eqcevstdx} \cite[Proposition 2.3.3 (1)]{das2009conditional},
and the standardization function $f$ satisfies \eqref{eqstdfn} by
Proposition \ref{propstdfn}.  

\subsubsection{Moment restrictions.}\label{subsubsec:restrict}
Section \ref{seccevmstd} considered  the standard case and found,
in particular, 
that if $X$ belongs to the standardized domain of attraction,  $0 \leq
\EP \xi \leq 1$ where $\xi\sim G$. 	 
When standardization is possible, a comparable moment restriction
occurs provided $X$ has a distribution in a domain of attraction.

Assume there exist normalizing functions $c(t)>0$ and $d(t) \in \mathbb{R}$
such that 
\begin{equation} \label{eqxdg}
t\P\bigg[\frac{X-d(t)}{c(t)} > x \bigg] \conv (1+\lambda x)^{-1/\lambda} \qquad x \in \E_\lambda, 
\end{equation} 
implying that $c,d \in \ERV_{\lambda,1}$ (see Section \refpc{secdoa}). 
If $(X,Y)$ follow a CEVM and \eqref{eqxdg} holds, then the vector
$(X,Y)$ belongs to a multivariate domain of attraction provided
$\lim_{t\goesto\infty} \alpha(t)/c(t) \in [0,\infty)$ \cite[Proposition 4.1]{das2011conditioning}.
Continuing the theme of assuming kernel convergence, 
consider the case  where \eqref{eqkernconvgen} holds under the same
normalization as in \eqref{eqxdg}:
\[ \mckernb{t}{[-\infty, c(t) x + d(t)]} \wc G(x).
\]
Then
from
 \eqrefpc{eqbterv}, 
 $d$ is a standardization function satisfying \eqref{eqstdfn}, and 
\[ \varphi(x) = \begin{cases} \inv{\lambda} ( x^\lambda - 1 ) &
  \lambda \neq 0 \\  \log x & \lambda = 0 \end{cases} .\] 
 Theorem \ref{thmcevstd} gives a standard CEVM for $(\ginv{d}(X),Y)$, and furthermore, $\ginv{d}(X)$ belongs to the standardized domain of attraction. Therefore, the distribution $G$ must satisfy 
\[ 
\int_0^\infty \P[\xi > \varphi(x)] dx  
 \leq 1. \]
Depending on $\lambda$, this reduces to  
\[ \begin{cases} \EP\xi^{1/\lambda} \ind{\xi>0} \leq \lambda^{-1/\lambda} & \lambda > 0  \\  
 \EP(-1/\xi)^{1/\abs{\lambda}} \ind{\xi<0} \leq \abs{\lambda}^{1/\abs{\lambda}} & \lambda < 0  \\ 
\EP e^\xi \leq 1	 &  \lambda = 0
\end{cases}. \]
Thus, we obtain a different condition for each class of extreme value distribution. 
In the Fr\'echet case, we have a bound on the $1/\lambda$-th moment of the right tail. 
If the domain of attraction is Weibull, this becomes an integrability condition near 0. 
Finally, in the Gumbel case, the right tail of $\xi$ is exponentially bounded, so all right-tail moments exist.


\subsection{Relation to the Heffernan and Tawn Model}

The CEVM of Theorem \ref{thmcevmgen} is inspired by the work of
Heffernan and Tawn \cite{heffernan2004conditional}. 
Where Heffernan and Tawn's model is based on the convergence of
conditional distributions as in \eqref{eqkernconvgen}, the general
CEVM
defined in \eqref{eqcevmdef}, \eqref{eqcondnondegennew} 
focuses on limits of joint distributions.  
Our Theorem \ref{thmcevmgen} shows that Heffernan and Tawn's
assumption \cite[Equation (3.1)]{heffernan2004conditional} leads to a
CEVM provided their normalization functions $\alpha$ and $\beta$ are 
ERV.   The fact that they require convergence \eqref{eqkernconvgen} 
 hold at all points $x$ suggests that they are expecting continuous
 limits whereas  we  framed the assumption as weak convergence. 

A condition such as
\eqref{eqkernconvgen}  tacitly assumes a
particular version of the conditional distribution. The issue of
version
cannot
be ignored, since Example \ref{exnotconv} shows that
\eqref{eqkernconvgen} holding for one particular version does not
imply that it holds for every version. 
The issue of version is usually handled by smoothness assumptions.

For a non-degenerate CEVM, the functions $\alpha$ and $\beta$
 are necessarily ERV. 
 Heffernan and Tawn assume a parametric form for these functions. 
They specify \[ \alpha(y) = b_{|i}(y) := y^{b_{|i}} = y^\rho \]
for some constant $\rho<1$ and
\[ \beta(y) = a_{|i}(y) := \begin{cases}
ay & \quad 0\leq\rho < 1,\ \ \text{with } a\in[0,1]\\
c-d\log y & \quad \rho<0\ \ \text{with } a=0,\, c\in\R,\, d\in [0,1]
\end{cases}. \]
Although more general models are possible, the form of the ERV limit
function $\psi$ in \eqrefp{eqervpsi} suggests that a parametric
approach is indeed reasonable.  


\section{General Normalizations for both $X$ and $Y$}
\label{seccevmgeny}

So far we  assumed that $Y$ satisfies  $t\P[Y>ty] \goesto \inv{y}$ for $y>0$. 
We now extend Theorem \ref{thmcevmgen} to the case
where $Y$ belongs to a general domain of attraction:  
\begin{equation} \label{eqydg}
 t\P[ Y > a(t) y + b(t) ] \conv (1+\gamma y)^{-1/\gamma} \qquad y\in \Eg, 
\end{equation}
and $\Eg := \{y: 1+\gamma y > 0\}$. 
Assume  $b(t)$ is given by \eqref{eqbtdf}. 
Without change,  \eqref{eqkernconvgen} may no longer be sufficient to obtain
a general CEVM limit \eqref{eqcevmdef}
if $Y$ requires normalization according to $a$ and $b$.

To relate kernel convergence to the CEVM when \eqref{eqydg} is the hypothesis,
there are two ways to proceed:  (i) Assume
$\mckernb{a(t) u + b(t)}{ [-\infty,\alpha(t)x + \beta(t)]} \goesto
\varphi_x(u)$ for $u>0$, and then \eqref{eqcevmdef} should follow from
arguments similar to those in Section \ref{seccevmgen}.  
(ii) Standardize $Y$ via the transformation $Y \mapsto \ginv{b}(Y)$,
use a version of $\P[X\in \cdot \given \ginv{b}(Y) = y]=:K^*(y,\cdot)$, and rely
on \eqref{eqkernconvgen} for $K^*$.
We show the consistency of these two approaches.


\subsection{Kernel Asymptotics}
\nc{\bi}{\ginv{b}}
\renewcommand{\bt}{\tilde{b}}
\nc{\bsi}{\ginv{{b^*}}}

The transition function $K :(-\infty,\infty) \times
\Bor[-\infty,\infty] \mapsto [0,1]$  is a specific
version of the conditional distribution of $X$ given $Y$, 
$ \mckern{y}{\cdot} = \P[X \in \cdot \given Y=y ]. $
To consider (ii)  above, we first express a version the conditional distribution of $X$ given $\bi(Y)$ in terms of $K$. 

When the distribution of $Y$ is not in the
Fr\'echet domain of attraction, the convergence \eqref{eqydg}, where $b$ is given
by \eqref{eqbtdf}, implies that  
$a,b \in \ERV_{\gamma,1}$ for some $\gamma \in \mathbb{R}$. Hence, $a\in\RV_\gamma$, and 
\begin{equation} \label{eqbtervr}
\frac{b(tx)-b(t)}{a(t)} \conv \begin{cases} 
\displaystyle\frac{x^\gamma - 1}{\gamma} & \gamma \neq 0 \\ \log x & \gamma = 0 
\end{cases}\,, \qquad x>0. 
\end{equation}
Inverting \eqref{eqbtervr} gives 
\begin{equation} \label{eqbtervinv}
\frac{\ginv{b}(a(t) x + b(t))}{t}\conv 
\begin{cases} (1+\gamma x)^{1/\gamma} & \gamma \neq 0 \\ e^x & \gamma = 0 \end{cases}\,, \qquad x \in \Eg. 
\end{equation}
Furthermore,  if $b^*$ is any function on $(0,\infty)$ satisfying 
\begin{equation} \label{eqbtequiv}
{{(b^*}(t)-b(t))}/{a(t)} \conv 0 \mt{as} t\goesto\infty, 
\end{equation}
then \eqref{eqydg}, \eqref{eqbtervr}, and \eqref{eqbtervinv} hold with $b$ replaced by $b^*$. 
A standard technique is to choose  a smooth, strictly monotone  $b^*$
 as is summarized next (cf. \cite{seneta1976regularly}).

\begin{lem} \label{lembtnice}
There exists a function $b^*$ satisfying \eqref{eqbtequiv} that is continuous and strictly monotone. 
\end{lem}

\begin{pf}
Consider cases:
If $\gamma=0$, then $b \in \Pi(a)$ and there exists 
\cite[Proposition 0.16]{resnick2007extreme}
$\bar{b}$
continuous, strictly increasing such that $(\bar{b}(t)-b(t))/a(t)
\goesto 1$. The choice
$b^*(x) = \bar{b}(\inv{e}x)$ satisfies \eqref{eqbtequiv}.  
If $\gamma>0$, then $b \in \RV_\gamma$, and $b(t)/a(t) \goesto \inv{\gamma}$ \cite[Theorem B.2.2 (1)]{de2006extreme}. 
Consequently, \mbox{\cite[Proposition 2.6 (vii)]{resnick2007heavy}} gives a continuous, strictly increasing function $b^* \sim b$. Writing  
\[ \frac{{b^*}(t)-b(t)}{a(t)} = \frac{b(t)}{a(t)} \bigg[\frac{b^*(t)}{b(t)} - 1 \bigg] \] shows that $b^*$ satisfies \eqref{eqbtequiv}. 
Finally, if $\gamma<0$, then $b(\infty) = \lim_{t\goesto\infty} b(t)$
exists finite, $b(\infty)-b \in \RV_\gamma$, and
$(b(\infty)-b(t))/a(t) \goesto -\inv{\gamma}$. Choose $\bar{b}$
continuous, strictly decreasing, with $\bar{b} \sim (b(\infty)-b)$,
and set $b^* = b(\infty)-\bar{b}$.  
\end{pf}

It is easier to deal with $b^*$ rather than $b$ since $\bsi(b^*(x)) =
b^*(\bsi(x)) = x$ but $b^*$ still standardizes $Y$.
By \eqref{eqbtervr}, $Y^* = \bsi(Y)$ is in the standard domain of
attraction when \eqref{eqydg} holds:  
\[ t\P[Y^* > ty] = t\P\bigg[\frac{Y-b^*(t)}{a(t)} > \frac{b^*(ty)-b^*(t)}{a(t)} \bigg] \conv \inv{y}, \qquad y>0. \]
Furthermore if $\mckernl{y}{\cdot} = \P[X\in \cdot \given Y=y]$, 
\begin{equation} \label{eqkerngenstd}
\mckern[K^*]{y}{\cdot} := \mckern{b^*(y)}{\cdot\,} 
\end{equation}
is a version of the conditional distribution $\P[X \in \cdot \given
Y^* = y]$. This follows from 
\begin{equation} \label{eqkernbt}
\P[X \in A , Y^* > y ] = \int_{(y,\infty)}  \mckern{b^*(u)}{A} \P[Y^*
\in du],
\quad \text{ measurable }A,\,y>0.
\end{equation} To see this write,
\ba
 \P[X \in A , &Y^* > y ] = \P[X \in A , Y > b^*(y)] = \int_{(b^*(y),\infty)}  \mckern{u}{A} \P[Y \in du] \\ 
 &=  \int_{(b^*(y),\infty)}  \mckernb{b^*(\bsi(u))}{A} \P[Y \in du],
\ea
using the fact that $b^*(\bsi(u))=u$ for all $u$. Finish with a change
variables to get \eqref{eqkernbt}. 

We now show that the two approaches to the CEVM discussed at the
beginning of Section \ref{seccevmgeny}, the direct approach and the
standardization approach, are  consistent. 

\begin{prop} \label{propkernunifgeny}
Suppose $Y$ has a distribution satisfying \eqref{eqydg} and  $K^*$ is given by \eqref{eqkerngenstd}. 
Given normalization functions $\alpha(t)>0$ and $\beta(t) \in \R$,
there exists a 
transition function $\phi^* : (0,\infty) \times \Bor[-\infty,\infty]
\mapsto [0,1]$   
such that, as $t\goesto\infty$, 
\begin{equation} \label{eqkernunifgenstd}
\mckernb[K^*]{tu_t}{[-\infty,\alpha(t) x+\beta(t)]}
\wc \mckern[\phi^*]{u}{[-\infty,x]} \mt{on} [-\infty,\infty]
\end{equation}
whenever $u_t \goesto u \in (0,\infty)$, iff
there is a 
transition function ${\phi : \Eg \times \Bor[-\infty,\infty] \mapsto [0,1]}$ 
such that, as $t\goesto\infty$, 
\begin{equation} \label{eqkernunifgeny}
 \mckernb{a(t)u_t+b(t)}{[-\infty,\alpha(t) x+\beta(t)]}  \wc \mckern[\phi]{u}{[-\infty,x]} \mt{on} [-\infty,\infty]
\end{equation}
whenever $u_t \goesto u \in \Eg$.
If these convergences hold, then 
\begin{enumerate}
\renewcommand{\labelenumi}{(\roman{enumi})}%
	\item $\alpha,\beta \in \ERV$; 
	\item $\phi^* = \tkg[G^*]$, a generalized tail kernel \eqref{eqtailkerngen} with $G^*=\mckernl[\phi^*]{1}{\cdot}$; 
	\item $\mckern[\phi]{u}{A} = \mckern[\tkg]{(1+\gamma u)^{1/\gamma}}{A}$, where $\tkg$ is a generalized tail kernel with $G=\mckernl[\phi]{0}{\cdot}$; and 
	\item the two transition functions are related by $G=G^*$. 
\end{enumerate}
\end{prop}

%
%
%
%
 %

\begin{pf}
Abbreviate $a_t = a(t)$ and $b_t = b(t)$. The convergences
\eqref{eqbtervr} and \eqref{eqbtervinv} are locally uniform on
$(0,\infty)$ (see 
Section \ref{secerv}).
Since $b^*$ satisfies \eqref{eqbtequiv}, it follows that 
\[ \frac{b^*(tu_t)-b_t}{a_t} \conv \frac{u^\gamma-1}{\gamma} \mt{whenever} u_t \goesto u \in (0,\infty), \]
and
\[\frac{\bsi(a_t u_t+b_t)}{t}\conv (1+\gamma u)^{1/\gamma} \mt{whenever} u_t \goesto u \in \Eg.  \] 
Assuming \eqref{eqkernunifgenstd}, for $u_t\goesto u \in \Eg$ we have 
\ba
\mckernb{a(t) u_t &+ b(t)}{[-\infty,\alpha(t) x+\beta(t)]} \\ 
 &= 
K\Big({b^*\big(t \{ \inv{t} \bsi(a_t u_t + b_t)\} \big)} \,,\, {[-\infty,\alpha(t) x+\beta(t)]}\Big) \\
 &= \mckernb[K^*]{t\{\inv{t} \bsi(a_t u_t + b_t ) \} } {[-\infty,\alpha(t) x+\beta(t)]} \\
 &\wc \mckernb[\phi^*]{(1+\gamma u)^{1/\gamma}}{[-\infty,x]} =: \mckern[\phi]{u}{[-\infty,x]}
\ea 
Conversely, if \eqref{eqkernunifgeny} holds, then for $u_t \goesto u>0$, 
\ba
\mckernb[K^*]{t &u_t}{[-\infty,\alpha(t) x+\beta(t)]} \\ 
 &= 
\mckernb{a_t \cdot \inv{a_t}({b^*(tu_t)- b_t} ) + b_t } {[-\infty,\alpha(t) x+\beta(t)]} \\
 &\wc \mckernb[\phi]{\inv{\gamma}(u^\gamma-1)}{[-\infty,x]} =: \mckern[\phi^*]{u}{[-\infty,x]}
\ea 
In either case, $G := \mckernl[\phi]{0}{\cdot} = \mckernl[\phi^*]{1}{\cdot} =: G^*$. Proposition \ref{propkernunifgen} shows that $\alpha$ and $\beta$ are ERV and $\phi^* = \tkg[G^*]$.  
Consequently, $\mckern[\phi]{u}{\cdot} = \mckern[\tkg]{(1+\gamma u)^{1/\gamma}}{\cdot}$. 
\end{pf}

Therefore, by Proposition \refp{propkernunifgen}, 
if there exists a non-degenerate distribution $G$ on $[-\infty,\infty)$ such that 
\begin{equation} \label{eqkernconvgenstd}
 \mckernb[K^*]{t}{[-\infty,\alpha(t) x + \beta(t)]} = \mckernb{b^*(t)}{[-\infty,\alpha(t) x + \beta(t)]} \wc G(x) 
\end{equation}
with $\alpha,\beta \in \ERV$, then \eqref{eqkernunifgeny} holds.  

How can we apply Proposition \ref{propkernunifgeny} starting from an assumption like \eqref{eqkernconvgenstd} on the kernel $K$ rather than $K^*$?
Because $b^*(\bsi(t)) = t$, \eqref{eqkernconvgenstd} can be written as 
\[ \mckernb{t}{[-\infty,\alpha\circ\bsi(t) x + \beta\circ\bsi(t)]} \wc
G(x) \mt{as} t \goesto y^*,  
\] 
where $y^*$ denotes the upper endpoint of the distribution of $Y$,
written as $y^* = \sup\{y : F_Y(y) < 1 \}$.  
Therefore, we require there to exist a non-degenerate distribution $G$ and normalization functions $\tilde{\alpha}>0$ and $\tilde{\beta}$ such that 
\begin{equation} \label{eqkernconvy}
 \mckernb{t}{[-\infty,\tilde{\alpha}(t) x+ \tilde{\beta}(t)]} \wc G(x) \mt{as} t \goesto y^*, 
\end{equation}
and $\alpha = \tilde{\alpha} \circ b^*$, $\beta = \tilde{\beta} \circ b^* \in \ERV$.

\subsection{CEVM Properties}

The standardization approach given in the previous section yields
 a CEVM when $Y$ belongs to a general domain of attraction.

\begin{thm} \label{thmcevmgeny}
Suppose $(X,Y)$ is a random vector on $\R^2$, where $F_Y \in
D(G_\gamma)$ according to \eqref{eqydg}
and 
$\mckernl{y}{\cdot} = \P[X \in \cdot \given Y=y]$
satisfies \eqref{eqkernconvy} for normalizing functions $\tilde{\alpha}>0$ and $\tilde{\beta}\in\R$ and non-degenerate limit distribution $G$ on $[-\infty,\infty)$. 
Let $b^*$ be the function satisfying \eqref{eqbtequiv} given by Lemma
\ref{lembtnice} and put $\alpha = \tilde{\alpha} \circ b^*$, $\beta =
\tilde{\beta} \circ b^*$.  
Then, as $t\goesto\infty$, 
\begin{equation} \label{eqcevgeny}
 t\P\bigg[\bigg(\frac{X-\beta(t)}{\alpha(t)}, \frac{Y-b(t)}{a(t)}\bigg)\in\cdot\,\bigg] \vconv \mu(\cdot) \mt{in} \mplus([-\infty,\infty]\times \Egc), 
\end{equation}
where $\mu$ is a non-null Radon measure satisfying the conditional non-degeneracy conditions \eqref{eqcondnondegennew},
iff $\alpha, \beta \in \ERV_{\rho,k}$.  
In this case, the limit measure $\mu$ is specified by 
\begin{equation} \label{eqlimmeasgeny}
\mu([-\infty,x]\times(y,\infty]) = 
\int_0^{(1+\gamma y)^{-1/\gamma}}   G( u^\rho x+\psi(u) )\,  du ,
\qquad x\in\R,\ y\in\Eg, 
\end{equation}
with $\psi$ as in \eqref{eqervpsi}. The expression \eqref{eqlimmeasgeny} is continuous in $x$ and $y$ if $(\rho,k) \neq (0,0)$. 
\end{thm}

\begin{pf}
First, observe that $Y^* = \bsi(Y) \in D(G^*_1)$.  
Defining the transition function $\mckernl[K^*]{y}{\cdot} = \P[X \in \cdot \given Y^*=y]$ as in \eqref{eqkerngenstd}, our hypotheses imply \eqref{eqkernconvgenstd}. 
Therefore, if $\alpha, \beta \in \ERV_{\rho,k}$, then  by Theorem \ref{thmcevmgen}, we have 
\begin{equation*} 
 t\P\bigg[\bigg(\frac{X-\beta(t)}{\alpha(t)}, \frac{Y^*}{t}\bigg)\in\cdot\,\bigg] \vconv \mu^*(\cdot) \mt{in} \mplus([-\infty,\infty]\times(0,\infty]), 
\end{equation*}
where $\mu^*$ is 
$\mu^*([-\infty,x]\times(y,\infty]) 
= \int_0^{\inv{y}} G(u^\rho x+\psi(u) )\,du  ,\,
x\in\R,\ y>0, $ and is
conditionally non-degenerate. 
Consequently, for $x\in\R$ and $y\in\Eg$, 
\ba 
t\P\bigg[\frac{X-\beta(t)}{\alpha(t)} 
\leq x,& \frac{Y-b(t)}{a(t)}> y\bigg] 
= t\P\bigg[\frac{X-\beta(t)}{\alpha(t)} \leq x, \frac{Y^*}{t} > \frac{\bsi(a(t)y+b(t))}{t}\bigg] \\ 
 &= \int_0^{(1+\gamma y)^{-1/\gamma}} G( u^\rho x+\psi(u) )\,du  = \mu([-\infty,x]\times(y,\infty]),  
\ea 
and the marginal transformation of $Y$ does not affect  conditional
non-degeneracy or continuity.  
Conversely, \eqref{eqcevgeny} implies that $\alpha,\beta\in \ERV$ \cite[Proposition 1]{heffernan2007limit}. 
\end{pf}

Alternatively, instead of standardizing $Y$, we could equally have used the convergence \eqref{eqkernunifgeny}, which holds under our assumptions by Propositions \ref{propkernunifgen} and \ref{propkernunifgeny}.

Recalling the forms of the limit measure given in Section \ref{seccevmgen}, we can express the limit measure in \eqref{eqlimmeasgeny} as 
\begin{align*} 
&\mu([-\infty,x]\times(y,\infty]) = \\
&\ \ \left\{\begin{alignedat}{2}
&\displaystyle\frac{1}{\rho \abs{x+k\inv{\rho}}^{1/\rho}}
\int_{0}^{\abs{x+k\inv{\rho}}(1+\gamma y)^{-\rho/\gamma}}
u^{(1-\rho)/\rho} G\bigl(u \sgn (x+k\inv{\rho}) &&- k\inv{\rho} \bigr) du  \\[-2mm]
&&& \phantom{\rho = 0,}\,\ \rho \neq 0  \\[2mm]
 &\displaystyle\frac{1}{\abs{k} e^{x/k}}\int_{-\infty}^{x\sgn(k) - \abs{k}\inv{\gamma}\log (1+\gamma y)} e^{u/\abs{k}} G( u\sgn(k)) du  
&& \rho = 0,\  
k \neq 0 \\[3mm]
&(1+\gamma y)^{-1/\gamma}G( x)
&& \rho = 0,\  
 k=0
\end{alignedat}\right.
\end{align*}
where $\sgn(v) = v/\abs{v} \ind{v \neq 0}$, and we 
read the measure as $(1+\gamma y)^{-1/\gamma} G(-k\inv{\rho})$ when $x=-k\inv{\rho}$ for the case $\rho \neq 0$. 

In Example \refp{excevmgen}, we presented a transition function satisfying \eqref{eqkernconvgen} which did not lead to a CEVM when paired with $Y\in D(G_1^*)$. 
We now show that a non-degenerate CEVM may be obtained if $Y$ belongs to a non-standardized domain of attraction.

\begin{exa}
Consider $Y\sim \text{Exp}(1)$, and $U \sim \text{Uniform}(0,1)$, independent of $Y$. Put $X = Ue^Y$. 
Note that $Y\in D(G_0)$ with $a(t)\equiv 1$, $b(t)=\log t$, since for $y\in\R$, 
\[ t\P(Y>y+\log t) = te^{-y-\log t} = e^{-y}. \] 
A version of the conditional distribution is given by 
\[ \mckern{y}{[0,x]} = \P[X\leq x \given Y=y] = \P[U\leq xe^{-y}] = xe^{-y} \smin 1. \]
Taking $\tilde{\alpha}(t) = e^t$, we saw in Example \ref{excevmgen} that
\[ \mckern{t}{\tilde{\alpha}(t)[0,x]} \wc x \smin 1 = G(x), \]
although $\tilde{\alpha}$ is not regularly varying. 
Since $b$ is continuous and strictly monotone, 
set $\alpha(t) = \tilde{\alpha}(b(t)) = t$. Then 
\[ \mckern[K^*]{t}{t[0,x]} = \mckern{b(t)}{\tilde{\alpha}(b(t))[0,x]} \wc G(x), \]
and $\alpha(t) \in \RV_1$. Hence, $\mckern[K^*]{tu}{\alpha(t)[0,x]}
\wc x\inv{u} \smin 1 = G(\inv{u}x)$, and $\mckernl[K^*]{y}{\cdot}
= \P[X \in \cdot \given e^Y = y]$.  
On the other hand, note that for $u\in \R$, 
\[ \mckernb{a(t)u + b(t)}{\alpha(t)[0,x]} 
= tx e^{-u-\log t} \smin 1 = xe^{-u} \smin 1 = G\big(\inv{(e^u)}[0,x]\big). \]
This illustrates the equivalence presented in Proposition \refp{propkernunifgeny}. 
Now, for $x>0$, $y > 0$, the joint distribution is given by 
\[ \P[X\leq x, Y>y] = \int_{\log x \smax y}^\infty x e^{-2u} du + \int_y^{\log x} e^{-u}du \ind{y<\log x}. 
\]
Therefore, for $x>0$, $y\in \R$, and large $t$, we have 
\[ t\P[X\leq tx, Y>y+\log t] = \left\{\begin{array}{lr}
\displaystyle\frac{xe^{-2y}}{2} & \text{if } \log x \leq y \\
\displaystyle  e^{-y} - \frac{1}{2x} & \text{if } \log x > y \\
\end{array}\right\} = \mu([0,x]\times(y,\infty]),  
\]
and $(X,Y)$ follow a CEVM by Theorem \ref{thmcevmgeny}. 
\end{exa}


\section{Conclusions and Future Directions}

In many statistical
contexts a conditional formulation such as \eqref{eqkernconvgen} is
convenient.  
An example is when we model $X$ as an explicit function of $Y$
or when we  work with distributions that have continuous densities,
in which case the natural choice of version of the conditional
distribution is the absolutely continuous one.
In such cases, Heffernan and Tawn \cite{heffernan2004conditional}
approach
is natural and  leads to a 
 parsimonious extremal model which can account for varying degrees of
 asymptotic independence.  
Heffernan and Tawn propose a semiparametric model, where the limit
distribution $G$ is estimated nonparametrically, and the normalization
functions $\alpha$ and $\beta$ belong to a parametric family.  
The extended regular variation of $\alpha$ and $\beta$ provides
justification for the form of the parametric family.
The formulas for the limit measure derived in our present paper
 show assuming conditional distributions leads to a simpler CEV model
 parametrized by the distribution $G$ and the pair $(\rho,k)$, along
 with $\gamma$, the extreme value index of $Y$.

Fitting a bivariate CEV model has been considered in
\cite{das2011detecting,fougeres2012estimation}.  
These authors discuss statistics for detecting a CEV model and
estimating the normalizing functions.  
However, open questions remain, such as the asymptotic
distributions of  estimators, and the appropriate method for
nonparametric estimation of $G$.

A natural extension of the bivariate model 
is to higher-dimensional
vectors and  this was the original intention of Heffernan and Tawn,
who apply their methodology to a five-dimensional air pollution
dataset.  
It is not clear how to
condition on more than one extreme variable; presumably there are connections
between such a model and the usual multivariate domain of attraction model.
Cases where asymptotic independence is present between
some pairs of variables but not others  requires careful
treatment and how to proceed in high dimensions is not apparent.

 \bibliographystyle{plain}
 \bibliography{condmod}

\end{document}